\documentclass[12pt]{article}
\usepackage{times}
\usepackage{a4wide}
\usepackage{amsfonts}
\usepackage{theorem}
\usepackage{enumerate}
\usepackage{amsmath}
\usepackage{amscd}
\usepackage{amssymb}

\theorembodyfont{\sl}
\newtheorem{conj}[subsection]{Conjecture}
\newtheorem{thm}[subsection]{Theorem}
\newtheorem{prop}[subsection]{Proposition}
\newtheorem{lem}[subsection]{Lemma}

\theorembodyfont{\rmfamily}









\newenvironment{prf}[1]{\trivlist
\item[\hskip \labelsep{\it
#1.\hspace*{.3em}}]}{~\hspace{\fill}~$\square$\endtrivlist}
\newenvironment{proof}{\begin{prf}{\bf Proof}}{\end{prf}}

\newcommand{\CC}{{\mathbb C}}
\newcommand{\RR}{{\mathbb R}}
\renewcommand{\SS}{{\mathbb S}}
\newcommand{\QQ}{{\mathbb Q}}
\newcommand{\FF}{{\mathbb F}}
\newcommand{\ZZ}{{\mathbb Z}}

\renewcommand{\AA}{{\mathbb A}}
\newcommand{\GG}{{\mathbb G}}

\newcommand{\Zhat}{\hat{\ZZ}}
\newcommand{\AAf}{\AA_{\rm f}}

\newcommand{\calX}{{\cal X}}

\newcommand{\bs}{\backslash}
\newcommand{\ti}{\times}
\newcommand{\ol}{\overline}

\newcommand{\lto}{\longrightarrow}

\newcommand{\Hom}{{\rm Hom}}
\newcommand{\Res}{{\rm Res}}
\newcommand{\Sh}{{\rm Sh}}
\newcommand{\Gal}{{\rm Gal}}
\newcommand{\ad}{{\rm ad}}

\newcommand{\discr}{{\rm disc}}

\newcommand{\inn}{{\rm inn}}
\newcommand{\GL}{{\rm GL}}

\newcommand{\MT}{{\rm MT}}

\newcommand{\Tr}{{\rm Tr}}
\newcommand{\Gm}[2]{{\GG_{{\rm m}{{#1}}}^{{#2}}}}

\begin{document}

\title{A conjecture of Yves Andr\'e. \footnote{The author was
supported by an EPSRC grant and the European Research Training Network
``Arithmetic Algebraic Geometry''. Address:
Imperial College, Dept. of Mathematics,
180~Queens Gate, SW7~2BZ London, UK. E-mail: andrei.yafaev@ic.ac.uk}}
\author{Andrei Yafaev} 
\date{\today}

\maketitle

\section{Introduction.}

In this article we deal with the following conjecture 
by Yves Andr\'e and Frans Oort.
    
\begin{conj}[Andr\'e-Oort]
Let $(G,X)$ be a Shimura datum. Let $K$ be a compact open subgroup
of~$G(\AAf)$ and let $S$ be a set of special points in
$\Sh_K(G,X)(\CC)$. Then every irreducible component of the Zariski
closure of $S$ in $\Sh_K(G,X)_\CC$ is a subvariety of Hodge type.
\end{conj}

The introduction to \cite{E1} (and references contained therein) contains a 
comprehensive exposition of terminology and notations relative to this 
conjecture.
Since we use the same terminology and notations, we do not reproduce them here.
The introduction to \cite{E1} also contains an exposition of results 
on this 
conjecture obtained before the article \cite{EY} came out.

In this article we prove the following theorem, which is actually the statement conjectured 
by Yves Andr\'e in 1989 in his book \cite{Andre} (see Problem 9). This statement
(without the assumption of the GRH) is now referred to
as a conjecture of Yves Andr\'e. 

\begin{thm} \label{main-thm}
Assume the Generalised Riemann Hypothesis (GRH) for CM fields.
Let $(G,X)$ be a Shimura datum. Let $K$ be a compact open subgroup of $G(\AAf)$.
Let $C$ be an irreducible closed algebraic curve contained in 
the Shimura variety $\Sh_K(G,X)$ and such that $C$ contains
an infinite set of special points.
Then $C$ is of Hodge type.
\end{thm}

In the article \cite{EY} we considered a curve in a Shimura variety 
$\Sh_K(G,X)$
containing an infinite set $S$ of special 
points satisfying the following condition.
There is a faithful rational representation of $G$ such that 
the $\QQ$-Hodge structures corresponding to the points in $S$ via this 
representation
lie in one isomorphism class. We proved that such a curve is of Hodge type.
This was the strongest result towards the Andr\'e-Oort conjecture at that time.

In the article \cite{EY} we introduced some technical tools 
to attack the Andr\'e-Oort conjecture.
In particular we obtained the following characterisation of subvarieties 
of Hodge type 
of a Shimura variety associated to a Shimura datum $(G,X)$ with $G$ semisimple 
of adjoint type.
Let  $Z$ be a Hodge generic subvariety of $\Sh_K(G,X)$ 
contained in its image by some Hecke correspondence $T_g$ with $g$ an element of 
$G(\QQ_p)$ i.e. $Z\subset T_gZ$. 
Suppose that $p$ is bigger than some integer depending on $G$, $X$, $K$ and $Z$ and
that $g$ is such that for any simple factor $G_i$ of $G$, the image of $g$ in
$G_i(\QQ_p)$ is not contained in a compact subgroup.
Then $Z$ is of Hodge type provided $Z$ contains at least one special 
point.

The strategy used to prove our main theorem \ref{main-thm} 
is the same as the 
one used 
in \cite{EY} (see Section 2 of \cite{EY} for details). 
We use the characterisation mentioned above.
After having reduced ourselves to the case where the group $G$ is 
semisimple of adjoint type and where the curve $C$ is Hodge generic,
we try to get $C$ to be contained in its image by a suitable
Hecke correspondence.
We consider intersections of $C$ with its images
$T_g C$ by Hecke correspondences $T_g$ with $g$ some elements of 
$G(\QQ_p)$ for various primes $p$. 
For suitably chosen $p$ and $g$ such intersection contains
a Galois orbit of some special point of $C$. 
We prove that one can choose a prime $p$ and an element $g$, both
satisfying the conditions mentioned above,
in such a way that the Galois orbit is too large for the intersection 
$T_gC \cap C$ to be finite.
The choice of a prime $p$ with this property is made possible by the 
assumption of the GRH 
and the use of the effective version of the Chebotarev density theorem.
We conclude that $C$ is of Hodge type.

The heart of this paper is a proof of a theorem about lower bounds for
Galois orbits of special points of Shimura varieties.
Our theorem on Galois orbits is a partial answer 
to Edixhoven's question Open Problem 14 in \cite{EMO}.
Using the GRH we refine lower bounds for Galois orbits given in \cite{EY}
enough to be able to prove the conjecture of Yves Andr\'e.

In section~2.2 we obtain precise information about Mumford-Tate groups of 
special points and their 
representations coming from special points on Shimura varieties.
This information allows us to bring the following improvement to the 
main result of \cite{EY}.

\begin{thm} \label{main-thm-2}
Let $(G,X)$ be a Shimura datum. Let $K$ be a compact open subgroup of $G(\AAf)$.
Let $C$ be an irreducible closed algebraic curve contained in 
the Shimura variety $\Sh_K(G,X)$ and such that $C$ contains
an infinite set $S$ of special points satisfying the following condition.

For any point $s$ of $S$ we choose an element $(\tilde{s},g)$ of 
$X\times G(\AAf)$ lying over $s$.
We suppose that 
the Mumford-Tate groups $\MT(\tilde{s})$ 
lie in one isomorphism class of $\QQ$-tori as $s$ ranges through the
set $S$.
Then $C$ is of Hodge type.
\end{thm}
 
\section{Lower bounds for Galois orbits.}

In this section we prove a theorem giving lower bounds
for Galois orbits of special points of Shimura varieties.

\begin{thm} \label{lower-bounds}
Assume the GRH for CM fields. Let $N$ be a positive integer.
Let $(G,X)$ be a Shimura datum with $G$ semi-simple of adjoint type,
and let $K$ be a neat compact open subgroup of 
$G(\AAf)$. Via a faithful representation of $G$, we view $G$ as a closed 
algebraic subgroup of $\GL_{n,\QQ}$, such that $K$ is contained in $\GL_n(\Zhat)$.
Let $V_\ZZ$ be the induced variation of $\ZZ$-Hodge structure on $\Sh_K(G,X)$.
For $s$ in $\Sh_K(G,X)$, we let $V_s$ be the corresponding Hodge structure and $\MT(V_s)$ its 
Mumford-Tate group (viewed as a closed algebraic subgroup of $\GL_{n,\ZZ}$).
Let $F\subset \CC$ be a number field over which $\Sh_K(G,X)$ admits a 
canonical model. For any special point $s$ in $\Sh_K(G,X)$, let $L_s$ be the splitting field 
of $\MT(V_s)$ and $d_{L_s}$ be the absolute value of its discriminant.

There exist real $c_1>0$ and $c_2>0$ such that for any special point $s$ in $\Sh_K(G,X)_F(\ol\QQ)$
we have:
$$
|\Gal(\ol\QQ/F){\cdot}s|>c_1 \log(d_{L_s})^N \prod_{\{\text{$p$
prime}\;|\;\text{$\MT(V_s)_{\FF_p}$ is not a torus}\}} c_2p 
$$
\end{thm}

\subsection{Reciprocity morphisms and Mumford-Tate groups.}

In this section we recall the definition of the Mumford-Tate group and 
reciprocity morphism attached to special elements of $X$ and prove some  
technical results about the Mumford-Tate groups and reciprocity morphisms 
to be used later on.

Let $(G,X)$ be a Shimura datum with $G$ semisimple of adjoint type,
let $V$ be a faithful rational representation of $G$
and let $V_\ZZ$ be a lattice in $V$.
Then $V_\ZZ$ induces a variation of $\ZZ$-Hodge structure over $X$.
Let $h$ be a special element of $X$. 
The morphism $h\colon \SS\lto G_\RR$, composed with the representation gives 
an $\RR$-Hodge structure $h\colon\SS\lto\GL(V_\RR)$.
Let $z$ and $\ol{z}$ be the generators of the character group of $\SS$.
The morphism $h$ corresponds to the decomposition
$$
V_\CC=\oplus_{p,q}V^{p,q}
$$
where $V^{p,q}$ is the $\CC$-vector subspace on which $\SS$ acts through the
character $z^p{\ol z}^q$. The spaces $V^{p,q}$ satisfy the following condition
 $\ol{V^{p,q}}=V^{q,p}$.
Let $W$ be the collection of pairs of integers $(p,q)$
that intervene in this representation.
Since the $\RR$-Hodge structures corresponding to elements of $X$ lie in one 
isomorphism class,
the set $W$ does not depend on the element $h$ in $X$.
The fact that $G$ is of adjoint type implies that for any $(p,q)$ in $W$ we have
$p+q=0$.
Let $M\subset\GL(V)$ be the Mumford-Tate group of $h$ and let $L$ be its splitting field.
Let us recall that $M$ is the smallest algebraic subgroup $H$ of $\GL(V)$ having the 
property that 
$h$ factors through $H_\RR$. The group $M$ is a $\QQ$-torus because $h$ is special 
and $M$ is given a $\ZZ$-structure by taking its Zariski closure in 
the $\ZZ$-group scheme $\GL(V_\ZZ)$.
We let $X^*(M)$ be the character group of $M$ i.e the group $\Hom(M_{\ol\QQ} , \Gm{\ol\QQ}{})$.
The group $X^*(M)$ is a free $\ZZ$-module of rank equal to the dimension of $M$ with
a continuous $\Gal(\ol\QQ / \QQ)$-action.

\begin{lem}
The field $L$ is a Galois CM field. 
Futhermore, the degree of $L$ is bounded in terms of the dimension of $V$.
\end{lem}

\begin{proof}
The field $L$ is Galois since it is the splitting field of a torus (the group $\Gal(\ol\QQ/L)$ is exactly the
kernel of the morphism $\Gal(\ol\QQ / \QQ)\lto \GL(X^*(M))$ hence is a normal subgroup of $\Gal(\ol\QQ/\QQ)$).
The fact that $L$ is a CM field follows from the fact that $\ad(h(i))$ is a Cartan involution of
$G_\RR$ (this is a part of axioms imposed upon Shimura data).

Let $E$ be the centre of the 
endomorphism algebra of the Hodge structure $V$.
The algebra $E$ is a finite product of number fields $E=E_1\ti\cdots \ti E_m$.
The torus $M$ is a subtorus of the torus $\prod_{i=1}^m \Res_{E_i/\QQ} \Gm{E_i}{}$. Hence $M$ is 
split over the composite of the Galois closures of the $E_i$ whose degree is 
clearly bounded in terms of 
the dimension of $V$ only.
\end{proof}

Let $T$ be the $\QQ$-torus $\Res_{L/\QQ}\Gm{L}{}$. 
Let $G_L$ be the Galois group of $L$ over $\QQ$ and
let $r\colon T \lto M\subset\GL(V)$ be the
reciprocity morphism associated to $h$.
Let us recall how $r$ is defined.
The morphism $h$ gives, by extending scalars from $\RR$ to $\CC$, 
the morphism $h_\CC$
from $\Gm{\CC}{}\ti \Gm{\CC}{}$ to $M_\CC \subset \GL(V_\CC)$.
Let $\mu\colon \Gm{\CC}{}\lto M_\CC$ be the morphism $h_\CC(z,1)$.
This morphism $\mu$ is defined over $L$.
Hence $\mu$ induces a morphism $\Gm{L}{}\lto M_L$, which, by taking the restriction of 
scalars from $L$ to 
$\QQ$ gives the morphism
$$
\Res_{L/\QQ}\mu\colon \Res_{L/\QQ}\Gm{L}{}\lto \Res_{L/\QQ}M_L.
$$
This morphism $\Res_{L/\QQ}\mu$ followed the by 
the norm morphism $\Res_{L/\QQ}M_L\lto M$ gives $r$.
The morphism $X^*(r)$ between character groups $X^*(M)$ and $X^*(T)$ 
is injective 
(because $r$ is a surjective morphism of $\QQ$-tori). 
The Galois module $X^*(T)$ is naturally isomorphic to $\ZZ[G_L]$.
We enumerate the elements of $G_L$ thus choosing a basis for $X^*(T)$ so
that it now makes sense to talk about coordinates of elements of $X^*(T)$.

The morphism $r$, when composed with the representation,  
defines an action of the torus $T$ on the $\QQ$-vector space $V$.
There is a subset $\calX$ of $X^*(T)$ such that this representation 
corresponds to a direct sum decomposition
$$
V_{\ol\QQ}=\oplus_{\chi\in\calX} V_{\chi}.
$$
where each $V_\chi$ is a $\ol{\QQ}$-subspace of $V_{\ol\QQ}$ 
on which $T_{\ol\QQ}$ acts through the character $\chi$.
The spaces $V_\chi$ satisfy the  condition that $V_\chi^\sigma=V_{\sigma\chi}$
(which insures that the representation is defined over $\QQ$).
  
The representation $h_\CC(z,1)$ of $\Gm{\CC}{}$ corresponds to the 
decomposition
$V_\CC=\oplus_{(p,q)\in W} V^{p,q}$
where $\Gm{\CC}{}$ acts via the character $z^p$ on $V^{p,q}$.
This representation is defined over $L$ hence induces a representation
of $\Gm{L}{}$. We get a decomposition $V_L=\oplus_{(p,q)\in W} V^{p,q}$
where $\Gm{L}{}$ acts through the character $z^p$ on $V^{p,q}$.
The representation $r$ of $T_L$ is obtained by
taking the restriction of scalars of this representation of $\Gm{L}{}$
followed by the norm from $L$ to $\QQ$.

It follows that the characters 
of $X^*(T)$ that belong to $\calX$ have coordinates 
(with respect to the basis we have chosen)
can be only integers 
$p$  or $q$ where 
${(p,q)}$
is some element of $W$.
In particular they are bounded, in absolute value, independently 
of the element $h$.
Furthermore, the characters in $\calX$ have the property that 
for any $\chi$ in $\calX$, the character $\chi \ol{\chi}$ is
the identity because 
the morphism $r$ satisfies the so-called Seere's condition (the group $G$ is of adjoint type)
and $p+q = 0$ for every pair $(p,q)$ in $W$. 
We refer to Section 2 of Chapter I of \cite{Milne}, in particular the Proposition 
2.4 for facts about Hodge structures of CM type.
We summarise what has been said in the following proposition.

\begin{prop} \label{coord}
There is an integer $k>0$ such that the following holds.
Let $h$ be a special element of $X$, $M$ its Mumford-Tate group 
and $L$ its splitting field.
Choose a basis for $X^*(T)$ by enumerating the elements of $G_L$.
With respect to the basis the coordinates of the 
characters of $T$ that
intervene in the decomposition 
$V_{\ol\QQ}=\oplus_{\chi\in\calX} V_{\chi}$
coming from the representation $r$ associated to $h$
have absolute value at most $k$.
Furthermore, for any character $\chi$ in $\calX$ the character 
$\chi\ol{\chi}$ is the identity.
\end{prop}

We now apply this Proposition to prove a number of results 
about Mumford-Tate groups of special elements of $X$ and reciprocity
morphisms attached to such elements.
These results will be used later on to prove lower bounds for Galois orbits.

\begin{prop} \label{index} 
There is a real $e>0$ such that the following holds.
Let $h$ be a special element of $X$.
Let $M$ be the Mumford-Tate group of $h$ and $L$ be its splitting field.
Let $r\colon T\lto M$ be the reciprocity morphism attached to $h$.
Let $p$ be a prime. 
The index of $r((\QQ_p\otimes L)^*)$ in $M(\QQ_p)$ is finite bounded  above 
by $e$.
The index 
of $r((\ZZ_p\otimes O_L)^*)$ in the maximal compact open subgroup
of $M(\QQ_p)$  is finite and bounded  above
by $e$.
\end{prop}
 
\begin{proof}
Let $P$ be the $\ZZ$-submodule of $X^*(T)$ spanned by
the vectors in $\calX$. Recall that we identify $X^*(M)$ with 
 its image by $X^*(r)$ i.e we view it as a submodule
of $X^*(T)$ and we have chosen a basis for $X^*(T)$. 
The module $X^*(M)$ is $P$.
The group $M(\QQ_p)$ is canonically isomorphic to
the group $\Hom_{G_L}(X^*(M),(\QQ_p\otimes L)^*)$
of $G_L$-invariant homomorphisms.
Similarly the group $(\QQ_p\otimes L)^*$ is
isomorphic to $\Hom_{G_L}(X^*(T),(\QQ_p\otimes L)^*)$
and the morphism $r\colon T(\QQ_p)\lto M(\QQ_p)$ is 
$$
r\colon  \Hom_{G_L}(X^*(T),(\QQ_p\otimes L)^*) \lto \Hom_{G_L}(X^*(M),(\QQ_p\otimes L)^*)
$$
which is just the restriction. 
The group $X^*(T)/P$ is a product 
of a free abelian group and a torsion group.
Let $E$ be the order of this torsion subgroup.
Since, by the previous Proposition, the coordinates of the vectors generating $P$ 
are bounded (in absolute value) by a unifrm constant $k$,
the number $E$ is bounded in terms of $k$ and $n_L$ only.
It is straightforward to see that
the order of the cokernel is bounded in terms of $E$ and $n_L$ 
only. The first claim follows.

The maximal compact open subgroup of $M(\QQ_p)$ is 
$\Hom_{G_L}(X^*(M),(\ZZ_p\otimes O_L)^*)$.
The second claim is proved using exactly the same arguments.
\end{proof}

\begin{prop} \label{coordinates}
There is an integer $B>0$ such that the following holds.
Let $h$ be a special element of $X$ and let $M$ be its Mumford-Tate group
and $L$ its splitting field.
Let $p$ be a prime splitting $L$ (hence $M$). There is a $\ZZ$-basis of the character group 
$X^*(M)$ such that the differences of coordinates
of the characters (with respect to this basis) that intervene in the representation $V_{\QQ_p}$ 
of $M_{\QQ_p}$ have absolute value at most $B$. 
\end{prop}

\begin{proof}
The module $X^*(M_{\QQ_p})$ is a submodule of $X^*(T)$ (along with its given basis).
generated by vectors whose coordinates are bounded in absolute value
by the integer $k$ from the Proposition \ref{coord}.
This integer is independent of the point $h$.
It follows that there is only finite number (depending on $k$ and $n_L$ only)
of possibilities for the set $\calX$ and hence for 
the submodule $X^*(M)$ of $X^*(T)$.
Choose some basis for $X^*(M)$ for each of this finite
number of cases. Take $B$ to be the maximum of absolute values of the 
differences of coordinates of characters in $\calX$ with respect to these bases.
\end{proof}
 
\begin{prop} \label{connected-comp}
There is a real $C>0$ such that the following holds.
Let $p$ be a prime.
For any special element $h$ in $X$ with Mumford-Tate group $M$ such that $M_{\FF_p}$ is a 
torus, the following holds.
Let $Y$ be a subspace of $V_{\FF_p}$. 
Let $T$ be the
stabiliser of $Y$ in $M_{\FF_p}$ (as defined in the Lemma 3.3.1 of \cite{EY}).
The order of the group of connected components of $T_{\ol\FF_p}$ 
has order at most $C$.
The order of the cokernel of the morphism $M(\FF_p)\lto (M/T)(\FF_p)$
is at most $C$.
\end{prop}
 
\begin{proof}
Proceeding as in the proof of the Lemma 4.4.1 of \cite{EY}, we reduce the proof of 
this proposition to the proof of the 
fact that stabilisers of lines satisfy the conclusion of the statement above.
We have a decomposition $V_{\ol{\FF_p}} = \oplus_{\chi\in \calX} V_\chi$.
Let $v$ be an element of $V_{\ol{\FF_p}}$, write
$v=\sum_{\chi} v_\chi$.
The stabiliser of the line $k v$ is the intersection of the 
kernels of $\chi-\chi'$ with $\chi$ and $\chi'$ 
distinct characters such that $v_\chi\not= 0$ and $v_{\chi'}\not= 0$.
Since the torsion of each $\ZZ$-module $X^*(T)/(\chi-\chi')\ZZ$ is 
bounded independently of $s$ and of the characters $\chi$ and $\chi'$ in
$\calX$, the order of the group of connected components 
of the stabiliser of $k \cdot v$  is bounded independently of 
$s$, $p$ and the subspace. This proves the first claim.

As for the second claim, using the Lemma 4.4.2 of \cite{EY}, we see that the 
order of the cokernel of the map $M(\FF_p)\lto (M/T)(\FF_p)$ is bounded
by the order of the group of connected components of $T_{\ol\FF_p}$ which is
uniformly bounded by what has just been said.
The second claim follows.
\end{proof}
 
\begin{prop} \label{roots}
There is a real $D>0$ such that the following holds.
Let $h$ be a special element of $X$ and let $M$ be its Mumford-Tate group.
Let $K_M$ be the maximal compact open subgroup of $M(\AAf)$.
The intersection $M(\QQ)\cap K_M$ is finite of order bounded by $D$.
\end{prop}

\begin{proof}
The group $M(\QQ)\cap K_M$ is finite because $M(\RR)$ is compact ($M(\RR)$ stabilises the
point $h$ of the Hermitian symmetric domain $X$ and the group $G_\RR$ is of adjoint type)
and the group $M(\QQ)\cap K_M$ is discrete. 
Let $L$ be the splitting field of $M$. 
Choose any basis for the character group $X^*(M)$ and use this basis to embed 
$M$ into a product of $\dim(M)$ copies of $T_L$.
Then the group $M(\QQ)\cap K_M$ is, via  this embedding, a finite subgroup of
the product of $d$ copies of $O_L^*$. 
It follows that it is contained in the
product of $\dim(M)$ 
copies of the group of roots of unity in $L$ which is finite of order bounded
independently of the point $h$. The claim follows.
\end{proof}
  
\subsection{Getting rid of $G$.}
  
Choose a set of
representatives $R$ in $G(\AAf)$ for the set of double classes
$G(\QQ)\backslash G(\AAf)/K$. Note that $R$ is finite. For $s$ in
$\Sh_K(G,X)$ there exists a unique $g_s$ in $R$ and an element
$\tilde{s}$ in $X$ unique up to $\Gamma_s:=G(\QQ)\cap g_sKg_s^{-1}$,
such that $s=\ol{(\tilde{s},g_s)}$.
Let $K_s$ be the compact open subgroup of $G(\AAf)$ defined by 
$K_s:=g_s K g_s^{-1}$.
We let 
$\MT(\tilde{s})$ be the 
Mumford-Tate group of $\tilde{s}$ (the smallest algebraic subgroup $H$ of $G$ such that $\tilde{s}$
factors through $H_\RR$).
The Mumford-Tate group $\MT(V_s)$ is the image of $\MT(\tilde{s})$ by the 
representation (this follows from the explicit construction of the 
variation of Hodge structures over $\Sh_K(G,X)$ given in the Section~3.2 
of \cite{EY}).
The element $\tilde{s}$ gives an embedding of the Shimura datum 
$(\MT(\tilde{s}), \{\tilde{s} \})$
into $(G,X)$. 

In this section we reduce the problem of giving a lower bound for the Galois orbit
of the point $s$ of $\Sh_K(G,X)$ to the one of giving 
a lower bound for the Galois orbit of the point $\ol{(\tilde{s},1)}$ of 
$\Sh_{K_s\cap \MT(\tilde{s})(\AAf)}(\MT(\tilde{s}))$.

\begin{prop}
The morphism of Shimura varieties
$$
\Sh_{K_s\cap \MT(\tilde{s})(\AAf)}(\MT(\tilde{s})) \lto \Sh_K(G,X)
$$ 
sending $\ol{(\tilde{s},t)}$ to $\ol{(\tilde{s},t\cdot g_s)}$
is injective.
\end{prop}
\begin{proof}
Let $M:=\MT(\tilde{s})$.
Let $H$ be the centraliser of $M$ in $G$.
Let $\ol{(\tilde{s},t)}$ and $\ol{(\tilde{s},t')}$ be two points of 
$\Sh_{K_s\cap \MT(\tilde{s})(\AAf)}(\MT(\tilde{s}))$
such that $\ol{(\tilde{s},t\cdot g_s)} = \ol{(\tilde{s},t'\cdot g_s)}$
in $\Sh_K(G,X)$. 
There exists an element $q$ of $H(\QQ)$ and an element $k$ of $K$ such that we have the following relation
$$
t = q t' g_s k g_s^{-1}
$$
Since $H(\QQ)$ and $M(\AAf)$ commute, this relation implies that
$tt^{'-1}$ belongs to $M(\AAf)\cap H(\QQ)U_s$ with $U_s:=H(\AAf)\cap K_s$.
Hence what we need to prove is that 
$M(\AAf)\cap H(\QQ)K_s = M(\QQ)(M(\AAf)\cap U_s)$. 
Consider the quotient of algebraic groups $H\lto \ol{H}=H/M$, which is well defined since $M$ 
is normal in $H$.
The image of $\ol{U_s}$ of $U_s$ in $\ol{H}(\AAf)$ is neat.
On the other hand $\ol{H}(\RR)$ is compact since $H(\RR)$ is compact (as a stabiliser of a point in
a hermitian symmetric domain and because $G_\RR$ is of adjoint type) and the map 
$H(\RR)\lto \ol{H}(\RR)$
is surjective on identity components. It follows that $\ol{H}(\QQ)$ is discrete in $\ol{H}(\AAf)$
and hence $\ol{H}(\QQ) \cap \ol{U_s}$ is trivial by neatness of $\ol{U_s}$.

Now suppose that $h$ is in $H(\QQ)$ and $u$ in $U_s$ such that $hu$ is
in~$M(\AAf)$. Then, in $\ol{H}(\AAf)$, we have
$\ol{h}{\cdot}\ol{u}=1$, hence $\ol{h}=\ol{u}=1$ in
$\ol{H}(\AAf)$. That means that $h$ is in $M(\QQ)$ and $u$ is in
$M(\AAf)\cap U_s$.
The claim follows.
\end{proof}

\subsection{Lower bounds for Galois orbits.}

We keep the notations of the preceding section.
Let furthermore $L$ be the splitting field of $\MT(\tilde{s})$.
Let $r$ be the reciprocity morphism attached to $\tilde{s}$ as explained in Section~2.2.
To simplify the notation we
write $M$ for $\MT(\tilde{s})$. 
The morphism $\Sh_{K_s\cap \MT(\tilde{s})(\AAf)}(\MT(\tilde{s})) \lto \Sh_K(G,X)$
is defined over $L$.
The action of $\Gal(\ol\QQ/L)$ on the Hecke orbit of
$(\tilde{s},g_s)$ is defined as follows.
The group $\Gal(\ol\QQ/L)$ acts through its maximal abelian quotient, which is,
by class field theory, isomorphic to a quotient of a product of a finite group 
of connected components of $(\RR\otimes L)^*$
(of order
bounded in terms of the degree of $L$ only) and of
$(\AAf\otimes L)^* / (\Zhat\otimes O_L)^*$.
Let $\sigma$ be an element of $\Gal(\ol\QQ/L)$ and
$t$ be an element of $(\AAf\otimes L)^*$ such that some element in the
 preimage of $\sigma$
in $(\AA\otimes L)^*$ followed by the projection to $(\AAf\otimes L)^*$
is $t$.
Then 
$$
\sigma\ol{(\tilde{s},g_s)} = \ol{(\tilde{s},r(t)\cdot g_s)}
$$
It follows that the size of the Galois orbit is, up to a uniformly
bounded factor, 
the size of the set $\ol{(\tilde{s}, r((\AAf\otimes L)^*)\cdot g_s)}$.
From the last lemma it follows that to prove the
Theorem \ref{lower-bounds}
it suffices to give a lower bound for
the size of the image of the set $(\tilde{s}, r((\AAf\otimes L)^*))$
in $\Sh_{K_s\cap M(\AAf)}(M)$.
Since the set $R$ of elements $g_s$
is finite, the index $[K:K_s\cap K]$ is bounded independently of $s$ and
it suffices to give a lower bound for the Galois orbit of the point $\ol{(\tilde{s},1)}$
of the Shimura variety $\Sh_{\GL_n(\Zhat)\cap M(\AAf)}(M)$.

\begin{lem} 
There is an element $q$ of $\GL_n(\QQ)$ such that the torus $M':=q M q^{-1}$ 
satisfies the condition that 
$M'_{\FF_p}$ is a torus for any prime $p$ not dividing the discriminant of 
$L$.
\end{lem}

\begin{proof}
Let $S$ be the finite set of primes $p$ such that $M(\ZZ_p)$ is not the
maximal compact subgroup $\Hom_{G_L}(X^*(M),(\ZZ_p\otimes O_L)^*)$ 
of $M(\QQ_p)$.
For every prime $p$ in $S$, choose a lattice $L_p$ in $\QQ_p^n$
invariant under the maximal compact subgroup of $M(\QQ_p)$.
Let $g=(g_p)$ be an element of $\GL_n(\AAf)$  
such that each
$g_p$ is an element of $\GL_n(\QQ_p)$ such that 
$L_p = g_p \ZZ_p^n$.
As $\GL_n(\AAf) = \GL_n(\QQ)\GL_n(\Zhat)$, we get an
 element $q$ of $\GL_n(\QQ)$ such that $g=q k$ for some $k$ in
 $\GL_n(\Zhat)$.
By the Lemma~3.3.1 of \cite{EY}, the torus $M':=q M q^{-1}$ is a torus for
every $p$ unramified in $L$.
\end{proof}
 
The morphism $\inn_q$ induces an isomorphism between
$M(\QQ) \bs M(\AAf) /M(\AAf)\cap \GL_n(\Zhat)$ and 
$M'(\QQ) \bs M'(\AAf) / M'(\AAf) \cap q\GL_n(\Zhat)q^{-1}$.
We let $r'$ denote the morphism $\inn_q \circ r$. 
To give a lower bound for the Galois orbit of the point $\ol{(\tilde{s},1)}$ of 
$\Sh_{ M'(\AAf) \cap q\GL_n(\Zhat)q^{-1}}(M')$ it suffices to
give a lower bound for the image of $r'((\AAf\otimes L)^*)$  in 
$M'(\QQ) \bs M'(\AAf) / M'(\AAf) \cap q\GL_n(\Zhat)q^{-1}$.

\begin{prop}
The size of the image of $r'((\AAf\otimes L)^*)$ in $\Sh_{M'(\AAf) \cap q\GL_n(\Zhat)q^{-1}}(M')$
is, up to a uniform (i.e depending only on the Shimura variety, not on $s$) constant, the size of the image of 
$r'((\AAf\otimes L)^*)$ in
$\Sh_{M'(\Zhat)}(M')$ times  
the size of the image of $r'((\AAf\otimes L)^*)\cap M'(\Zhat)$ in
$M'(\Zhat)/M'(\Zhat)\cap q \GL_n(\Zhat) q^{-1}$.
\end{prop}
 
\begin{proof}
We are interested in the size of the set 
$$r'((\AAf\otimes L)^*) / r'((\AAf\otimes L)^*) \cap (M'(\QQ)(q\GL_n(\Zhat)q^{-1}\cap M'(\AAf))).$$
Since $M'(\Zhat)$ is the maximal compact subgroup of $M'(\AAf)$, we have an inclusion 
$$
M'(\AAf)\cap q\GL_n(\Zhat)q^{-1}\subset M'(\Zhat).
$$
Hence the size of the set we are interested in is the size of 
$$r'((\AAf\otimes L)^*) / r'((\AAf\otimes L)^*)\cap M'(\QQ)M'(\Zhat)$$
times that of 
$$
r'((\AAf\otimes L)^*)\cap M'(\QQ)M'(\Zhat) / 
r'((\AAf\otimes L)^*)\cap (M'(\QQ) M'(\AAf) \cap  q\GL_n(\Zhat)q^{-1}).
$$
The order of  $M'(\QQ)\cap M'(\Zhat)$ is
 bounded independently of the point $s$ by the Proposition \ref{roots} hence the size of 
$$r'((\AAf\otimes L)^*)\cap M'(\QQ)M'(\Zhat) / r'((\AAf\otimes L)^*) \cap 
(M'(\QQ)q\GL_n(\Zhat)q^{-1}\cap M'(\AAf))$$
is, up to a uniformly bounded constant, that of 
the of the image of $r'((\AAf\otimes L)^*)\cap M'(\Zhat)$ in
$M'(\Zhat)/M'(\Zhat)\cap q \GL_n(\Zhat) q^{-1}$. 
\end{proof}

\begin{prop}
The size of the image of $r'((\AAf\otimes L)^*)\cap M'(\Zhat)$ in
$M'(\Zhat) / M'(\Zhat) \cap q\GL_n(\Zhat)q^{-1}$
is at least
$$
\prod_{\{\text{$p$ prime}\;|\;\text{$\MT(\tilde{s})_{\FF_p}$ is not a torus}\}} c_2p 
$$
where  $c_2$ is strictly positive real constant independent of $p$ and $s$.
\end{prop}

\begin{proof}
Clearly it suffices to give a lower bound for the size of the image
of $r'((\Zhat\otimes O_L)^*)$ in
$M'(\Zhat) / M'(\Zhat) \cap q\GL_n(\Zhat)q^{-1}$.

The group $M'(\Zhat) / M'(\Zhat) \cap q\GL_n(\Zhat)q^{-1}$ is a product over primes $p$
such that $M'(\ZZ_p) \not= M'(\ZZ_p) \cap q\GL_n(\ZZ_p)q^{-1}$ of groups
$M'(\ZZ_p) / M'(\ZZ_p) \cap q\GL_n(\ZZ_p)q^{-1}$.
Hence we need for every such $p$ to give a lower bound for the image of $r'((\ZZ_p\otimes O_L)^*)$ in
the quotient. 

Let $p$ be such prime.
Using the fact that $r'((\ZZ_p\otimes O_L)^*)$ has uniformy bounded index in 
$M'(\ZZ_p)$, we see that up to a uniformly bounded constant, the size of the 
image of  $r'((\ZZ_p\otimes O_L)^*)$ in $M'(\ZZ_p) / M'(\ZZ_p) \cap q\GL_n(\ZZ_p)q^{-1}$
is the size of $M'(\ZZ_p) / M'(\ZZ_p) \cap q\GL_n(\ZZ_p)q^{-1}$.

The size of the set $M'(\ZZ_p) / M'(\ZZ_p) \cap q\GL_n(\ZZ_p)q^{-1}$
is the size of the orbit $M'(\ZZ_p) \cdot q_p \ZZ_p^n$
where $q_p$ is the image of $q$ in $\GL_n(\QQ_p)$. 
The Proposition 4.3.9 of \cite{EY} tells us that the size of this orbit is at least $c p$ if 
$M'_{\ZZ_p}$ does not fix the lattice $q_p \ZZ_p^n$. This last condition means exactly that
$\MT(V_s)_{\FF_p}$ is not a torus (by Lemma 3.3.1 of \cite{EY}).
The only thing we have to do is to check that the constant $c$ can actually be taken 
independent of $s$. The proof of the Proposition~4.3.9 of \cite{EY}
tells us that this constant depends only on the 
orders of the groups of connected components of the stabilisers of
subspaces of $\FF_p^n$ in $M_{\FF_p}$.
By the Lemma \ref{connected-comp},
these orders are bounded independently of $p$ and $s$ and the subspace in
question.
\end{proof}

The remaining task in order to prove the theorem \ref{lower-bounds}
is to give a lower bound for the size of the 
set $r' ((\AAf\otimes L)^*) / r' ((\AAf\otimes L)^*)\cap M'(\QQ)M'(\Zhat)$.
We prove the following Theorem.

\begin{thm}
Assume the GRH for CM fields. Let $N$ be a positive integer.
There is a real constant $c>0$ independent of the choice of $s$ and $M'$ (but depending on $N$)
such that the size of the set
 $r'((\AAf\otimes L)^*) / r'((\AAf\otimes L)^*)\cap M'(\QQ)M'(\Zhat)$
 is at least $c \log(d_L)^N$.
\end{thm}

\begin{proof}
In what follows, we write $M$ for $M'$ and $r$ for $r'$ to simplify the notations.
Let $n_L$ be the degree of $L$ over $\QQ$.
Let $m>0$ be an integer at most $\frac{\log(d_L)^5}{15n_L \log\log(d_L)}$ and let 
$p_1,\dots,p_m$ be $m$ distinct primes split in $L$ and 
smaller than $\log(d_L)^5$. 
Their existence
is provided by the effective Chebotarev theorem (under GRH), provided $d_L$ is bigger than
some absolute constant, which we assume.
We refer to the Proposition 8.2 of \cite{E1} for the exact statement of the
effective Chebotarev theorem that we use.
For each $i=1,\dots,m$, we choose a place $v_i$ of $L$ lying over $p_i$.
We let $P_i$ be the uniformiser at the place $v_i$.
Let $n_1,\dots, n_m$ be integers satisfying $|n_i|<N$.
Let $I$ be the
element of $(\AAf\otimes L)^*$ that equals $P_i^{n_i}$ at the place $v_i$ for 
$i=1,\dots,m$ and $1$ elsewhere. Suppose that $r(I)$ belongs to $M(\QQ)M(\Zhat)$.
Let $\pi$ be a corresponding element of $M(\QQ)$ (this element is defined up to
an element of $M(\QQ)\cap M(\Zhat)$ which is, by the Proposition \ref{roots}, 
a finite group of uniformly bounded order). Let, as before, 
$T$ be the torus $\Res_{L/\QQ}\Gm{L}{}$.
Recall that there is a set $\calX\subset X^*(T)$ such that
the representation $r$ of $T$ gives a decomposition 
$V_L=\oplus_{\chi\in\calX} V_\chi$ which is Galois invariant.
Let $d$ be the cardinality of $\calX$ and let 
$(\pi_1,\dots,\pi_d)$ be the $d$ elements of $L^*$ which are images of $\pi$ by the $d$ characters
in $\calX$. The field $\QQ(\pi_1,\dots,\pi_d)$ is Galois because the set $\calX$ is Galois invariant.

\begin{lem}
Suppose that not all $n_i$ are zero.
Then the field $\QQ(\pi_1,\dots,\pi_d)$ is $L$.
\end{lem}

\begin{proof}
It suffices to prove that the group $\Gal(L/\QQ(\pi_1,\dots,\pi_d))$ acts 
trivially on $\QQ\otimes X^*(M)$
(alternatively on $\QQ\otimes X_*(M)$, $X_*(M)$ being the group of cocharacters).

To simplify the exposition we suppose that $m=1$ 
(the general case is done using exactly the same arguments).
Let $\sigma$ be an element of $\Gal(L/\QQ(\pi_1,\dots,\pi_d))$.
So we have a prime $p$ splitting $L$, we choose  
a place $v$ of $L$ lying over $p$ and a uniformiser $P$ at $v$.
We consider the idele $I=P^n$ with $n>1$ some integer. We suppose that $r(I)$ belongs to $M(\QQ)M(\Zhat)$.
As $p$ splits $M$, we have 
$M(\QQ_p)=\Hom(X^*(M),\QQ_p^*)=X_*(M)\otimes \QQ_p^*$.
It follows that the evaluation map $v_p\colon \QQ_p^*\lto\ZZ$  induces an
isomorphism between 
 $M(\QQ_p)/M(\ZZ_p)$ and the group of cocharacters $X_*(M)$ of $M$.
Let $K$ be the kernel of $r$, then we have an exact sequence of $\QQ$-vector
spaces with 
$\Gal(\ol\QQ/\QQ)$-action
$$
0\lto \QQ\otimes X_*(K)\lto \QQ\otimes X_*(T) \lto \QQ\otimes X_*(M) \lto 0.
$$
It suffices to prove that $\sigma$ acts trivially on 
$\QQ\otimes (X_*(T)/X_*(K))$.
Since $\sigma$ fixes each $\pi_i$ and the set of $\pi_i$ is $G_L$-invariant,
$\sigma$ fixes all the elements $r(\tau I)$ of $M(\QQ_p)/M(\ZZ_p)$
with $\tau$ ranging through $G_L$.
The Galois action on $M(\QQ_p)/M(\ZZ_p)$ is being given by 
identifying it with $X_*(M)$ which has a Galois action.
Since the morphism $X_*(T)\lto X_*(M)$ is surjective, 
for any $\tau$ in $G_L$ we have
$\sigma \tau I = \tau I$ in  $\QQ\otimes (X_*(T)/X_*(K))$.
Let $e_1,\dots , e_{n_L}$ be the basis of $\QQ\otimes X_*(T)$ given by the
the $n$th powers of uniformisers at the places lying over $p$.
Their images in $\QQ\otimes (X_*(T)/X_*(K))$ generate this vector space.
Since $\sigma$ fixes these elements,
 $\sigma$ acts trivially on $\QQ\otimes (X_*(T)/X_*(K))$.
The claim follows.
\end{proof}

Let $x$ be the integer $(p_1\cdots p_m)^{Nk}$ with $k$ the integer from the Proposition \ref{coord}.
Let $\chi$ be a character in $\calX$.
The element $x \chi(I)$ of $(\AAf\otimes L)^*$ belongs to $\Zhat \otimes O_L$.
On the other hand this element is of the form $x \pi_i$ (for some $i$) times some element of
$(\Zhat\otimes O_L)^*$. It follows that $x \pi_i$ is in $O_L$.
We replace $\pi_i$ with $x \pi_i$.
The fact that $\chi\ol\chi$ is the identity shows that 
$\pi_i\ol{\pi_i}$ is $x^2=(p_1\cdots p_m)^{2Nk}$.
The field $\QQ(\pi_1,\dots,\pi_d)$ is $L$.
Let us choose a basis $b_1,\dots, b_{n_L}$ of $L$ over $\QQ$ consisting of monomials
in $\pi_1,\dots, \pi_d$ of degree bounded by a constant depending on $n_L$ only.
The discriminant of the ring $\ZZ[b_1,\dots,b_{n_L}]$ is the discriminant
of the matrix $A$ whose entries $A_{ij}$ are $A_{ij}=\Tr_{L/\QQ}(b_i b_j)$.
The absolute values of the $A_{ij}$ are bounded by a uniform power of $(p_1\cdots p_m)^N$.
We see that the discrimiant of $A$ is the sum of $n_L !$ terms 
whose absolute values are bounded by a uniform power of $(p_1\cdots p_m)^N$
hence there is a uniform constant $t$ such that 
$$
|\discr \ZZ[b_1,\dots , b_{n_L}]|\leq  (p_1\cdots p_m)^{Nt}.
$$
Since $\ZZ[b_1,\dots,b_{n_L}]$ is an order in $O_L$, we have
$$
|\discr\ZZ[b_1,\dots,b_{n_L}]|\geq d_L
$$
Replacing $t$ with $5t$, we get
the following inequality
$$ 
\log(d_L)^{Nmt} > d_L
$$
Hence, if $I$ is such that  $r(I)$ belongs to $M(\QQ)M(\Zhat)$, then 
$$
Nm >\frac{\log(d_L)}{t\log\log(d_L)}
$$
Let us now consider elements 
of $(\AAf\otimes L)^*$ that equal $P_i^{n_i}$ at the place $v_i$ and
$1$ outside of the places $v_i$ and 
where 
$|n_i|<N/2$ with $N$ and $m$ are such that 
$Nm \leq \frac{\log(d_L)}{t\log\log(d_L)}$.
From the above inequality it follows that these elements have distinct non-trivial images 
in $M(\AAf)/M(\QQ)M(\Zhat)$ by $r$.
It follows that the set $r((\AAf\otimes L)^*) / r((\AAf\otimes L)^*)\cap M(\QQ)M(\Zhat)$
contains at least $(N/2)^m$ elements if $N$ and $m$ being such that 
$Nm \leq \frac{\log(d_L)}{t\log\log(d_L)}$.   
Taking $m=\frac{\log(d_L)}{20tNn_L\log\log(d_L)}$ (which is possible by the effective Chebotarev),
we easily see that $(N/2)^m$ is at least $c \log(d_L)^N$ elements where $c$ is 
some real positive constant not depending on $s$ (but of course depending on $N$).
\end{proof} 

\section{Proof of main results.}
In this section we prove the Theorems \ref{main-thm} and \ref{main-thm-2}.

Let $(G,X)$ be a Shimura datum and $K$ a compact open subgroup in $G(\AAf)$.
Let $C$ be an irreducible closed algebraic curve in $\Sh_K(G,X)$ containing 
an infinite set $\Sigma$ of special points.
For any special point $s$ of $C$ we let $L_s$ be the splitting field of 
the Mumford-Tate group of some element $\tilde{s}$ lying over $s$ 
and we let $d_s$ be the absolute value of the discriminant of $L_s$.
 
\begin{prop}
Suppose that the discriminant of $L_s$ is bounded as $s$ ranges through $\Sigma$.
Then $C$ is of Hodge type.
\end{prop}

\begin{proof}
We can assume that $G$ is semisimple of adjoint type (passing to the adjoint group does not
change the property of $C$ being of Hodge type by the Proposition~2.2 of \cite{EY} and 
does not change the property that $d_{L_s}$ is bounded).
Let us choose some faithful representation $V$ of $G$.
Since the discriminant of $L_s$ is bounded as $s$ ranges through $\Sigma$,
 there are only finitely many possibilities for $L_s$.
Hence we can assume that for all points $s$ in $\Sigma$, the field $L_s$ is the same field $L$.
For any $s$ in $\Sigma$, we let $\tilde{s}$ be an element of $X$ such that $s=\ol{(\tilde{s},g)}$
for some $g$ in $G(\AAf)$.
The reciprocity morphism $r_{\tilde{s}}$ gives a rational representation of the torus $T:=\Res_{L/\QQ}\Gm{L}{}$.
This representation corresponds to a direct sum decomposition $V_L=\oplus_{\chi\in \cal{X}} V_\chi$ 
for some subset $\cal{X}$ of $X^*(T)$.
As before, we identify the $G_L$-module $X^*(T)$ with $\ZZ[G_L]$ and enumerate elements of $G_L$ thus 
getting a basis for $X^*(T)$.
Using the fact that coordinates of the characters in the set $\calX$
are bounded in absolute value by $k$ which does not depend on $s$ (Proposition \ref{coord}), we 
see that there are only finitely many possibilities for the
set $\cal{X}$ as $s$ ranges through $\Sigma$.
Hence, possibly replacing $\Sigma$ by an infinite subset, we can and do assume that the set $\calX$
 is constant as $s$ ranges through 
$\Sigma$. We can further assume that the dimensions of the $V_\chi$ are constant.
We now see that
the $\QQ$-Hodge stuctures $V_{\tilde{s}}$ are isomorphic 
as $s$ ranges through $\Sigma$.
Hence $C$ is of Hodge type by the main theorem of \cite{EY}.
\end{proof}

From the proof of this Proposition the Theorem \ref{main-thm-2} follows. 
Indeed, let $C$ be a curve in $\Sh_K(G,X)$
that contains an infinite set $\Sigma$ of points such that the corresponding
Mumford-Tate groups are isomorphic as $\QQ$-tori.
Since the Mumford-Tate groups 
of points of $\Sigma$ are isomorphic, they have the same splitting field.
From the proof of the above proposition, it follows that $C$ contains an infinite set 
of special points such that the $\QQ$-Hodge structures 
corresponding to these points via some faithful representation of $G$
lie in one isomorphism class. By the main result of \cite{EY},
$C$ is of Hodge type.

In what follows we assume that $d_{L_s}$ is 
unbounded as $s$ ranges through $\Sigma$.
From Propositions 2.1 and 2.2 of \cite{EY}, it follows that we can assume $G$ to be semisimple of adjoint type 
and $C$ to be Hodge generic.
Write $G=G_1\ti\cdots\ti G_r$ where $G_i$ are simple. We can and do assume that $K$ is the product
of compact open subgroups $K_p$ of the $G(\QQ_p)$ and that $K$ is neat.
Choose a faithful representation $V$ of $G$ through which we view $G$ as a closed subgroup of $\GL_{n,\QQ}$
such that $K$ is in $\GL_n(\Zhat)$.
Also choose a $K$-invariant lattice in $V_{\AAf}$.
This gives a variation of $\ZZ$-Hodge structure on $\Sh_K(G,X)$ 
(Section 3.2 of \cite{EY}).
Let $X^+$ be a connected component of $X$. After possibly having replaced $C$ by an irreducible component
of its image under a suitable Hecke correspondence we can and do assume that $C$ is contained in the image 
$S$ of $X^+\ti\{ 1 \}$ in $\Sh_K(G,X)$.
Since $C$ contains an infinite set of special points which are in $\Sh_K(G,X)(\ol\QQ)$, $C$ is defined
over a Galois number field $F$ containing the reflex field of $(G,X)$ 
(as an absolutely irreducible closed subscheme 
$Z_F$ of
$\Sh_K(G,X)_F$). 

\begin{prop}
Assume the GRH for CM fields.
There is
a prime $p$ and a point $s$ in $\Sigma$ which have the following properties
\begin{enumerate}
\item $p$ splits $\MT(V_s)$.
\item $\MT(V_s)_{\FF_p}$ is a torus.
\item Let $k$ be an integer as in the Corrollary 7.4.4 of \cite{EY}.
Then $| \Gal(\ol\QQ/F)\cdot s|>p^k$.
\end{enumerate}
\end{prop}
\begin{proof}
Let, as in the section 7 of \cite{EY}, define the function $i\colon\Sigma\lto \ZZ$ as follows.
For $s$ in $\Sigma$, let $i(s)$ be the number of prime numbers $p$
such that $\MT(V_s)_{\FF_p}$ is not a torus.
Then, by the Theorem \ref{lower-bounds},
there exist real $c_1>0$ and $c_2>0$ such that for any $s$ in $\Sigma$ we have 
$$
|\Gal(\ol\QQ/F)\cdot s|>c_1 \log(d_{L_s})^{5k} c_2^{i(s)}i(s)!
$$
where $k$ is the integer from the Corollary 7.4.4 of \cite{EY}.
Using this inequality and the effective Chebotarev theorem (in the form stated in the 
Proposition 8.2 of \cite{E1}) we see that
the number of primes
split in $L_s$ and smaller than $|\Gal(\ol\QQ/F)\cdot s|^{1/k}$ is bigger than $i(s)$ when $d_{L_s}$ is 
large enough. This finishes the proof of the proposition.
\end{proof}
Take a prime $p$ and a point $s$ given by the previous proposition.
Let $m$ be an element of $G(\QQ_p)$ given by the Corollary 7.4.4 of \cite{EY}
(this Corrolary can be applied because of our Proposition \ref{coordinates}).
Then some Galois conjugate of $s$ is in $C\cap T_mC$ and since $C\cap T_mC$ is defined 
over $F$ the whole Galois orbit of $s$ is contained in $C\cap T_mC$.
If the intersection $C\cap T_mC$ was finite, its cardinality would be smaller than $p^k$.
By the choice of $p$ and $s$, this intersection can not be finite hence $C$ is contained in $T_mC$.
We conclude that $C$ is of Hodge type using the Theorem 7.1 of \cite{EY}.

\section*{Acknowledgements}
This work is a continuation of author's Phd thesis carried out under supervision of Bas Edixhoven.
The author is grateful to him for many important discussions on the subject
and many importants comments on previous versions of this paper.
The author is grateful to Ralph Greenberg who provided him with his unpublished notes
on abelian varieties with complex multiplication. 
The author is grateful to Alexei Skorobogatov for some useful discussions.
The author is grateful to the EPSRC and Arithmetic Algebraic Geometry network for
financial support.

\end{document}